\documentclass[a4paper,11pt]{article}

\usepackage{amsmath}
\usepackage{amsfonts}
\usepackage{amssymb}
\usepackage{amsthm}
\usepackage{array}
\usepackage[utf8]{inputenc}
\usepackage{setspace}

\newtheorem{Theo}{Theorem}

\newtheorem{Conj}{Conjecture}

\title{Proof of a Combinatorial Conjecture Coming from the
PAC-Bayesian Machine Learning Theory.}
\author{Malik Younsi}
\begin{document}

\maketitle
\begin{abstract}
\setlength \parindent{0cm}
\noindent
We give a proof of a conjecture of A. Lacasse in his doctoral thesis \textbf{\cite{LAC}} which has applications in machine learning algorithms. The proof relies on some interesting binomial sums identities introduced by Abel $(1839)$, and on their generalization to the multinomial case by Hurwitz $(1902)$.
\end{abstract}

\section{The conjecture}
In his thesis \textbf{\cite{LAC}}, A. Lacasse gives the following conjecture :

\begin{Conj}
For $m \in \mathbb{N}$, define
$$\xi(m):= \sum_{k=0}^m \binom{m}{k} \left( \frac{k}{m} \right)^k \left( 1-\frac{k}{m} \right)^{m-k}$$
and
$$\xi_2(m):= \sum_{j=0}^m \sum_{k=0}^{m-j} \binom{m}{j} \binom{m-j}{k} \left( \frac{j}{m} \right)^j \left( \frac{k}{m} \right)^k \left( 1-\frac{j}{m}-\frac{k}{m} \right)^{m-j-k}.$$
Then
$$\xi_2(m) = m+\xi(m) \qquad (m \in \mathbb{N}).$$
\end{Conj}
This conjecture has applications in Machine Learning Theory, see \textbf{\cite{LAC}}. It was verified numerically for $m$ up to $4000$.

\section{Proof of the conjecture}

To prove the conjecture, we first rewrite the functions $\xi(m)$ and $\xi_2(m)$ under a more convenient form.

Define
\begin{eqnarray*}
\alpha(m) := m^m \xi(m) &=& \sum_{k=0}^m \binom{m}{k} m^k \left( \frac{k}{m} \right)^k m^{m-k}  \left( 1-\frac{k}{m} \right)^{m-k}\\
&=& \sum_{k=0}^m \binom{m}{k} k^k (m-k)^{m-k}
\end{eqnarray*}
and

\begin{eqnarray*}
\beta(m) &:=& m^m \xi_2(m)\\
 &=& \sum_{j=0}^m \sum_{k=0}^{m-j} \binom{m}{j} \binom{m-j}{k} m^j\left( \frac{j}{m} \right)^j m^k\left( \frac{k}{m} \right)^k m^{m-j-k}\left( 1-\frac{j}{m}-\frac{k}{m} \right)^{m-j-k}\\
&=& \sum_{j=0}^m \sum_{k=0}^{m-j} \binom{m}{j} \binom{m-j}{k} j^j k^k (m-j-k)^{m-j-k}.
\end{eqnarray*}

$\alpha(m)$ and $\beta(m)$ are sums of binomial and multinomial type, respectively. The conjecture is thus equivalent to the following :
$$\beta(m)-\alpha(m) = m^{m+1} \qquad (m \in \mathbb{N}).$$

Some numerical experimentations (including consultation of the On-line Encyclopedia of Integer Sequences) seem to suggest the following identities :

\begin{equation}
\label{eq1}
\alpha(m) = \sum_{j=0}^m m^j \frac{m!}{j!} \qquad (m \in \mathbb{N}),
\end{equation}

\begin{equation}
\label{eq2}
\beta(m) = \sum_{j=0}^m m^{m-j} \binom{m}{j} (j+1)! \qquad (m \in \mathbb{N}).
\end{equation}

Note that if $\textbf{(\ref{eq1})}$ and $\textbf{(\ref{eq2})}$ hold, then the conjecture holds, as can be seen by an elementary calculation :
\begin{eqnarray*}
\beta(m) - \alpha(m) &=& \sum_{j=0}^m m^{m-j} \binom{m}{j} (j+1)! - \sum_{j=0}^m m^j \frac{m!}{j!}\\
&=& \sum_{k=0}^m m^k  \binom{m}{m-k} (m-k+1)! - \sum_{j=0}^m m^j \frac{m!}{j!}\\
&=& \sum_{k=0}^m m^k  \frac{m!}{k!}  (m-k+1) - \sum_{j=0}^m m^j \frac{m!}{j!}\\
&=& \sum_{k=0}^m m^k \frac{m!}{k!} (m-k)\\
&=& m\sum_{k=0}^m m^k \frac{m!}{k!} - \sum_{k=0}^m km^k \frac{m!}{k!}\\
&=& \sum_{j=1}^{m+1} m^j \frac{m!}{(j-1)!} - \sum_{k=1}^m m^k \frac{m!}{(k-1)!}\\
&=& m^{m+1}
\end{eqnarray*}

After some research in the literature of combinatorial identities, we found identities $\textbf{(\ref{eq1})}$ and $\textbf{(\ref{eq2})}$ (under a slightly different form) in \textbf{\cite{RIO}}.
\\

More precisely, consider $\textbf{(\ref{eq1})}$. Define, for $m \in \mathbb{N}$, $x,y \in \mathbb{R}$, $p,q \in \mathbb{Z}$ :
$$A_m(x,y;p,q):= \sum_{k=0} ^m \binom{m}{k} (x+k)^{k+p}(y+m-k)^{m-k+q}.$$
The case $p=-1,q=0$ is well known : it is the so called \textit{Abel's Binomial Theorem}. Our case of interest is $x=0,y=0,p=0,q=0$. In \textbf{\cite{RIO}}, p.$21$, we find the identity

$$A_m(x,y;0,0) = \sum_{k=0}^m \binom{m}{k} k! (x+y+m)^{m-k}.$$
With $x=0,y=0$, this gives
$$\alpha(m)=A_m(0,0;0,0) = \sum_{k=0}^m \binom{m}{k} k! m^{m-k} = \sum_{j=0}^m \frac{m!}{j!} m^j,$$
which is the required identity $\textbf{(\ref{eq1})}$.
\\

For identity $\textbf{(\ref{eq2})}$, we need a multinomial version of $\textbf{(\ref{eq1})}$. This can be found in  \textbf{\cite{RIO}}, p.$25$, equation $(35)$.

For $x_1, x_2, \dots x_n \in \mathbb{R}$ and $p_1, \dots, p_n \in \mathbb{Z}$, define
$$A_m(x_1, \dots, x_n; p_1, \dots, p_n):= \sum \frac{m!}{k_1 ! k_2! \dots! k_n !} \prod_{j=1}^n (x_j+k_j)^{k_j+p_j},$$
where the sum is taken over all integers $k_1, \dots, k_n$ with $k_1 + \dots + k_n = m$.
\\

Clearly, our case of interest is $n=3$, $x_1=x_2=x_3 =0$, $p_1=p_2=p_3=0$ :

$$\beta(m) = A_m(0,0,0;0,0,0) \qquad (m \in \mathbb{N}).$$

Equation $(35)$,p.$25$ of \textbf{\cite{RIO}} gives the formula
$$A_m(x_1, \dots, x_n; 0, \dots, 0) = \sum_{k=0}^m \binom{m}{k} (x_1+x_2+\dots+x_n+m)^{m-k} \alpha_k(n-1),$$
where $\alpha_k(r):= \frac{(r+k-1)!}{(r-1)!}$.

For $n=3$, we have that $\alpha_k(n-1) = (k+1)!$ and with $x_1=x_2=x_3=0$, the above formula becomes
$$\beta(m) = A_m(0,0,0;0,0,0) = \sum_{k=0}^m \binom{m}{k} m^{m-k} (k+1)!,$$
which is exactly $\textbf{(\ref{eq2})}$.

We summarize all this in the following Theorem :

\begin{Theo}
For $m \in \mathbb{N}$, define
$$\xi(m):= \sum_{k=0}^m \binom{m}{k} \left( \frac{k}{m} \right)^k \left( 1-\frac{k}{m} \right)^{m-k}$$
and
$$\xi_2(m):= \sum_{j=0}^m \sum_{k=0}^{m-j} \binom{m}{j} \binom{m-j}{k} \left( \frac{j}{m} \right)^j \left( \frac{k}{m} \right)^k \left( 1-\frac{j}{m}-\frac{k}{m} \right)^{m-j-k}.$$

Then we have
$$\xi(m) = \frac{1}{m^m} \sum_{j=0}^m m^j \frac{m!}{j!} \qquad (m \in \mathbb{N})$$
and
$$\xi_2(m) = \frac{1}{m^m}\sum_{j=0}^m m^{m-j} \binom{m}{j} (j+1)! \qquad (m \in \mathbb{N}).$$

Furthermore,
$$\xi_2(m) = \xi(m)+m \qquad (m \in \mathbb{N}).$$
\end{Theo}

To conclude, we want to emphasize on the fact that the above not only proves the conjecture but also gives simpler expressions for the functions $\xi(m)$ and $\xi_2(m)$. These expressions are more convenient to handle numerically.

\bibliographystyle{amsplain}

\end{document}